\documentclass[11pt]{article}
\usepackage{bez123,calc,curves,ebezier,epic,eepic,graphicx,multiply,rotating}
\everymath{\displaystyle}
\usepackage{graphicx}
\usepackage{graphicx}
\usepackage{subfigure}
\usepackage{pst-all}
\usepackage{color}
\usepackage{amssymb}
\usepackage{amsmath}
\newtheorem{prethm}{{\bf Theorem}}

\newenvironment{thm}{\begin{prethm}{\hspace{-0.5
               em}{\bf .}}}{\end{prethm}}
\newtheorem{prelemma}{{\bf Lemma}}

\newenvironment{lemma}{\begin{prelemma}{\hspace{-0.5
               em}{\bf .}}}{\end{prelemma}}
\newtheorem{preex}{{\bf Example}}

\newtheorem{preprop}{{\bf Proposition}}

\newtheorem{precor}{{\bf Corollary}}

\newenvironment{cor}{\begin{precor}{\hspace{-0.5
               em}{\bf .}}}{\end{precor}}
\newtheorem{preremark}{{\bf Remark}}

\newtheorem{preprob}{{\bf Problem}}

\newtheorem{predefin}{{\bf Definition}}

\newtheorem{preconj}{{\bf Conjecture}}

\newtheorem{preprobb}{{\bf Problem}}

\newtheorem{prelem}{{\bf Theorem}}

\newenvironment{proof}{{\bf Proof.}\rm }{\hfill{$\Box$}}

\newtheorem{presolution}{{\bf Solution.}}

\def\newpic#1{}
\def\qed{\ifhmode\unskip\nobreak\fi\quad\ifmmode\Box\else$\Box$\fi}
\title{\vspace{-0.1cm}\Large\bf On monopoly and dynamic monopoly of Cartesian product of graphs with constant thresholds}
\author{\large\bf Nahideh Asadi~~~~~~Manouchehr Zaker\footnote{E-mail: mzaker@iasbs.ac.ir}
\vspace{5mm}\\
    Department of Mathematics,\\
     Institute for Advanced Studies in Basic Sciences,\\
     Zanjan 45137-66731, Iran}

    \date{}

\begin{document}
\maketitle\
\begin{abstract}
\noindent Let $G$ be any simple and undirected graph. By a threshold assignment $\tau$ in $G$ we mean any function $\tau:V(G)\rightarrow \mathbb{N}$ such that $\tau(v)\leq d_G(v)$ for any vertex $v$ of $G$. Given a graph $G$ with a threshold assignment $\tau$, a subset of vertices $M$ is said to be a $\tau$-monopoly if there exist at least $\tau(v)$ neighbors in $M$ for any vertex $v \in V(G) \setminus M$. Similarly, a subset of vertices $D$ is said to be a $\tau$-dynamic monopoly if starting with the set $D$ and iteratively adding to the current set further vertices $u$ that have at least $\tau(u)$ neighbors in it, results in the entire vertex set $V(G)$. Denote by $mon_{\tau}(G)$ (resp. $dyn_{\tau}(G)$) the smallest cardinality of a $\tau$-monopoly (resp. $\tau$-dynamic monopoly) of the graph among all others. In this paper we obtain some lower and upper bounds for these two parameters with constant threshold assignments for Cartesian product graphs. Our bounds improve the previous known bounds. We also determine the exact value of these two parameters with fixed thresholds in some Cartesian graph products including cycles and complete graphs.
\end{abstract}

\noindent {\bf Mathematics Subject Classification}: 05C69, 05C76, 05C35.

\noindent {\bf Keywords:} Dominating sets; Monopolies; Dynamic monopolies; Cartesian product.

%%%%%%%%%%%%%%%%%%%%%%%%%%%%%%%%%%%%%%%%%%%%%%%%%%%%%%%%%%%%%%%%%%%%%%%%%%%

\section{Introduction}

\noindent All graphs in this paper are simple and undirected. Let $G$ be a graph with the edge set $E(G)$. If two arbitrary vertices $v, v'$ are adjacent in $G$ then we write $vv'\in E(G)$. We refer the reader to the book \cite{W} for general notations and concepts of graph theory. In this paper we deal with Cartesian product of graph. We begin with its definition. The Cartesian product of $G$ and $H$, written $G\Box H$, is the graph with vertex set $V(G)\times V(H)$ specified by putting $(u, v)$ adjacent to $(u', v')$ if and only if (l) $u = u'$ and $vv'\in E(H)$, or (2) $v = v'$ and $uu'\in E(G)$. The vertex set of $G\Box H$ can be displayed as a $|V(G)| \times |V(H)|$ array of vertices, in which the rows (resp. columns) are indexed by the vertices of $G$ (resp. $V(H)$). The subgraph of $G\Box H$ induced by each row (resp. column) is isomorphic to the graph $H$ (resp. $G$). We will use this representation of the vertices of $G\Box H$ throughout of paper. For more notations and knowledge concerning Cartesian product of graphs we refer the reader to the book \cite{IK}. We assume that the vertices of $G$ are labeled by natural numbers according to a threshold assignment $\tau: V(G)\longrightarrow \mathbb{N}$, where for each vertex $\tau(v)\leq d_G(v)$, where $d_G(v)$ stands for the degree of $v$ in $G$. A subset of vertices $M$ is said to be a $\tau$-static monopoly ($\tau$-monopoly, for simplicity) if every vertex $v\in V(G)\setminus M$ has at least $\tau(v)$ neighbors in $M$. We denote the size of smallest $\tau$-monopoly set by $mon_{\tau}(G)$ which is called the $\tau$-monopoly number of $G$. Monopoly of graphs have been widely studied in the literature, e.g. \cite{FKRRS, KNSZ, KPY2, MR, SZ}. A threshold assignment $\tau$ is called simple majority (resp. strict majority) if for every vertex $\tau(v)=d_G(v)/2$ (resp. $\tau(v)=(d_G(v)+1)/2$), where $d_G(v)$ is the degree of $v$ in $G$. Monopoly of graphs with respect to the strict majority assignment was also studied by some authors under the terminology of global offensive alliance. We refer the reader to \cite{FFGHHKLS} and also to \cite{YR} for offensive alliances in Cartesian product graphs. The paper \cite{KPY1} provides results on monopolies for direct product of graphs. The other concept similar to $\tau$-monopoly when $\tau$ is a constant function is $k$-domination sets. We refer to \cite{S} for results on $k$-domination number of some Cartesian product graphs.

\noindent A subset $ D $ of $V(G)$ is said to be a $\tau$-dynamic monopoly (or simply a $\tau$-dynamo) if the vertices of $ G $ can be partitioned into subsets $ D_{0},D_{1}, \ldots,D_{k} $ such that $ D_{0} = D $ and for any $ i = 1, \ldots, k - 1  $ each vertex $ v\in D_{i+1} $ has at least $ \tau(v) $ neighbors in $ D_{0}\cup \cdots \cup D_{i} $. The dynamic monopolies formulate the spread of influence in social networks, where a social network is formulated by a graph in which the participants and social links are represented by vertices and edges of the graph, respectively (see e.g. \cite{DR, FKRRS}). In fact the threshold of each vertex is interpreted as the level of susceptibility of the vertex with respect to a certain influence in the social network. Because of this application of dynamic monopolies in phenomena of spread of influence, the vertices of $D_0$ are called initially activated vertices. For each $v\in D_i$, we say $v$ is activated at $i$th time step. By the $\tau$-dynamic monopoly number we mean the cardinality of the smallest $\tau$-dynamic monopoly set of $G$. By a $dyn_\tau(G)$-set $D$ we mean any dynamic monopoly for $G$ of cardinality $dyn_\tau(G)$.
Dynamic monopolies have been widely investigated in the literature \cite{ABW, DR, FKRRS, KSZ, Z}. Also dynamic monopolies of graph products was studied in \cite{ABST, ATZ, CHY}. If in a graph $G$, all vertices have the same threshold, i.e. for some $t$ and for all $v\in V(G)$, $\tau(v)=t$, then we say that $G$ has a constant threshold $t$. In this case we simply write $mon_t(G)$ (resp. $dyn_t(G)$) instead of $mon_{\tau}(G)$ (resp. $dyn_{\tau}(G)$). The study of monopolies and dynamic monopolies with respect to simple or strict majority is very common in this area. Throughout the paper for any graph $G$ and any $S\subseteq V(G)$, we denote by $G[S]$ the subgraph of $G$ induced by the vertices of $S$.

\noindent {\bf The paper is organized as follows}. Section 2 devotes to monopolies. In this section we determine the exact value of the smallest monopoly with arbitrary but fixed thresholds for some Cartesian graph products including cycles and complete graphs. A lower bound is also given for a general case. In Section 3 we determine the smallest dynamic monopoly in some Cartesian product graphs and also obtain some lower and upper bounds for this parameter, which improve the previous known bounds.

\section{Monopolies in Cartesian graph products}

\noindent The first result concerns the toroidal graphs i.e. the Cartesian product of two cycles. When the underlying graph is 4-regular (such as toroidal graphs) then the simple majority threshold is equivalent to the threshold assignment in which all thresholds are 2. For this reason we obtain some bounds for such threshold assignments. In the following theorem we denote by $v_{ij}$, $1\leq i,j \leq n$, the vertex located in the $i$th horizontal copy and the $j$th vertical copy of $C_n \Box C_n$.

\begin{thm}
For the monopoly number of $C_n \Box C_n$ with a constant threshold $2$ we have
\begin{center}
$mon_2 (C_n \Box C_n) \left \{
\begin{array}{lll}
= \dfrac{n^2}{3} \qquad \quad \ \quad \quad \quad  \quad \quad n=3t  & \\
\\
\leq \dfrac{(n-1)(n+3)}{3} \quad\  \quad \quad n=3t+1  & \\
\\
\leq \dfrac{(n+1)^2}{3}-1 \ \quad \quad \quad \quad n=3t+2.
\end{array}
\right.$
\end{center}
\end{thm}

\noindent \begin{proof}
We have three cases. If $n=3t$ select the vertices
\begin{flushleft}
$M=\{v_{13}, v_{22}, v_{31}\}\cup \{v_{16}, v_{25}, \ldots, v_{61}\}\cup \ldots \cup \{v_{1(n-3)}, v_{2(n-4)}, \ldots, v_{(n-3)1}\} \cup \{v_{1n}, v_{2(n-1)}, \ldots, v_{n1}\}\cup \{v_{4n}, v_{5(n-1)}, \ldots, v_{n4}\}\cup \{v_{7n}, v_{7(n-1)}, \ldots, v_{n7}\}\cup \ldots \cup \{v_{(n-2)n}, v_{(n-1)(n-1)}, v_{n(n-2)}\}.$
\end{flushleft}
It is easily seen that each vertex in $(C_n \Box C_n)\setminus M$ is dominated by two vertices of $M$. The number of selected vertices in this case is
\begin{center}
$n+2[(n-3)+(n-6)+\cdots +3]=\dfrac{n^2}{3}.$
\end{center}
The $2$-static monopoly number of the graph in this case is exactly $n^2/3$ as every vertex of $G \setminus M$ is dominated by only $2$ monopoly vertices and also each monopoly vertex contributes in activating all of its four neighbors.\\
Now, if $n=3t+1$, then select
\begin{flushleft}
$M'=\{v_{13}, v_{22}, v_{31}\}\cup \{v_{16}, v_{25}, \ldots, v_{61}\}\cup \ldots\cup \{v_{1(n-1)}, v_{2(n-2)}, \ldots, v_{(n-1)1}\}\cup \{v_{3n}, v_{4(n-1)}, \ldots, v_{n3} \cup \{v_{4(n-1)}, v_{5(n-2)}, \ldots, v_{(n-1)4}\}\cup \{v_{7(n-1)}, v_{7(n-2)}, \ldots, v_{(n-1)7}\}\cup \ldots \cup \{v_{(n-1)n}, v_{n(n-1)}\}.$
\end{flushleft}
Set $M'$ along with $\{v_{1n}, v_{4n}, \ldots, v_{(n-3)n}\}\cup \{v_{n1}, v_{n4}, \ldots, v_{n(n-3)}\}$
is a $2$-static monopoly with the total number of vertices
\begin{center}
$(2+3+5+6+8+9+\cdots + n-1)+2 \dfrac{(n-1)}{3}=\dfrac{(n-1)(n+3)}{3}.$
\end{center}
If $n=3t+2$, then we apply a slight different method for selecting the vertices. The monopoly vertices are as follow.
\begin{flushleft}
$M''=\{v_{12}, v_{21}\} \cup \{v_{15}, v_{24}, \ldots, v_{51}\}\cup \{v_{18}, v_{27}, \ldots, v_{81}\}\cup \{v_{1n}, v_{2(n-1)}, \ldots, v_{n1}\}\cup \{v_{4n}, v_{5(n-1)}, \ldots, v_{n4}\}\cup \ldots \cup \{v_{(n-1)n}, v_{n(n-1)}\}$
\end{flushleft}
along with $\{v_{3n}, v_{6n}, \ldots, v_{(n-2)n}\}\cup \{v_{n3}, v_{n6}, \ldots, v_{n(n-2)}\}$
represents a $2$-static monopoly. The number of these vertices is
\begin{center}
$2(2+5+8+\cdots + n-3)+n+ 2\dfrac{(n-2)}{3}=2\dfrac{(n-2)}{6}(n-1)+n+2\dfrac{(n-2)}{3}=\dfrac{(n+1)^2}{3}-1.$
\end{center}
\end{proof}

\noindent The next theorem concerns the monopoly number of $C_m\Box K_n$.

\begin{thm}
For any cycle graph $C_m$ and any complete graph $K_n$, \footnotesize{\begin{equation*}
	mon_t (C_m \Box K_n)\left \{ \begin{array}{lll} = m(t-2) \qquad \quad \quad \quad \quad t-1\leq n\leq 2(t-2),\ m\ is\ even & \\ \leq (m-3)(t-2)+2n \quad \ t-1\leq n\leq 2(t-2),\ m\ is\ odd & \\ = \dfrac{mn}{2} \quad \quad \quad \quad \qquad \quad \quad 2(t-2)+1\leq n\leq 2t,\ m\ is\ even & \\ = \dfrac{(m-1)n}{2}+t-1 \quad \quad \ 2(t-2)+1\leq n\leq 2t,\ m\ is\ odd & \\ \leq mt \ \quad \qquad \qquad \qquad \quad\ n> 2t. \end{array} \right.
	\end{equation*}}
\end{thm}

\noindent \begin{proof} We first make a monopoly of the desired size. There are three cases:\\ \textbf{1)} $t-1\leq n\leq 2(t-2)$. If $m$ is an even integer we select $t-2$ vertices from the left side of the first row and the same number from the right side of the next row and continue this way of selection for the next odd and even rows. The total number of monopoly vertices in this case is $m(t-2)$.\\ Now, if $m$ is an odd then we select $t-2$ vertices in the same manner as above from rows $1$ to $m-1$. From the last row we choose $n-(t-2)$ vertices on its left side and the same number of vertices from its right. The number of monopoly vertices in this case is $(m-1)(t-2)+2n-2(t-2)=(m-3)(t-2)+2n$.\\ \\ \textbf{2)} $2(t-2)+1\leq n\leq 2t$. In this case from $i$th row select $t-i$ vertices on the left and $n-t+i$ vertices from the next row on its right-hand ($i$ being the smallest) in a way that $t-i$ and $n-t+i$ have a difference of at most $1$. If $\vert t-i-(n-t+i) \vert=\vert t-i-n+t-i \vert= \vert 2t-2i-n \vert =1$ then these two numbers might be either $(a)$ $t, t-1$ or $(b)$ $t-1, t-2$ respectively. There will be $(mn)/2$ monopoly vertices when $m$ is even. If $m$ is an odd integer we do the same process except for the last row for which if case $(a)$ happens we select $t-1$ vertices from the right side. On the other hand, if case $(b)$ happens we select $t-1$ vertices from the left. The number of monopoly vertices in each of the two cases is \begin{center} $ \dfrac{(m-1)n}{2}+t-1.$ \end{center} See Figure \textbf{1} for this case. Note that when $n$ is odd, $t-i$ and $n-t+i$ can not be equal. Otherwise, $t-i=n-t+i$ which leads to $n=2(t-i)$, a contradiction.\\ \\ \textbf{3)} Now if $n\geq 2t$, then we choose $t$ vertices from each row. In other words, there will be $mt$ monopoly vertices.\\ \\ In case $(1)$ and $(2)$ we prove the necessity of the theorem.
	\textbf{1):}
	It is necessary to select at least $t$ vertices from each row. Assume to the contrary that there exists a row, $l$, in the product graph from which at most $t-3$ vertices belong to the monopoly set. Considering that any vertex in the graph is adjacent to exactly two vertices in its column, one can verify that any vertex $v \in l$ is adjacent to at most $t-1$ active vertices and this is a contradiction.\\
	\textbf{2):}
	Assume that there are $x$ and $m-x$ rows with $t-2$ and $n-t+2$ initially activated vertices, respectively. Hence, \begin{center} $x(t-2) + (m - x)(n - t + 2) = x(2t - 4 - n) + m(n - t + 2)$ \end{center} Now, let $m$ be an odd integer, then the above statement \begin{center} $ \geq \dfrac{m-1}{2}(2t -4 - n) + mn - mt + 2m =\dfrac{m+1}{2}(n) - t + 2\geq \dfrac{m-1}{2}(n) +t-1$ \end{center} The last equality comes from the fact that $n\geq 2(t-2)+1$.\\ In case that $m$ is even, we have the statement above to be \begin{center} $\geq \dfrac{m}{2}(2t - 4 - n) + mn - mt + 2m =\dfrac{mn}{2}.$ \end{center}
\end{proof}

\begin{figure}[ht]
\begin{center}
{
\textcolor{red}{\begin{tabular}
{|p{.25cm}|p{.25cm}|p{.25cm}|
p{.25cm}|p{.25cm}|p{.25cm}|p{.25cm}|p{.25cm}|p{.25cm}|p{.25cm}|}
\hline{\bf $*$}&{\bf $*$}&{\bf $*$}&{\bf $*$}&{\bf $*$}&&&&\\[+0.2eM]
\hline&&&&&{\bf $*$}&{\bf $*$}&{\bf $*$}&{\bf $*$}\\[+0.2eM]
\hline{\bf $*$}&{\bf $*$}&{\bf $*$}&{\bf $*$}&{\bf $*$}&&&&\\[+0.2eM]
\hline&&&&&{\bf $*$}&{\bf $*$}&{\bf $*$}&{\bf $*$}\\[+0.2eM]
\hline{\bf $*$}&{\bf $*$}&{\bf $*$}&{\bf $*$}&{\bf $*$}&&&&\\[+0.2eM]
\hline&&&&&{\bf $*$}&{\bf $*$}&{\bf $*$}&{\bf $*$}\\[+0.2eM]
\hline{\bf $*$}&{\bf $*$}&{\bf $*$}&{\bf $*$}&{\bf $*$}&&&&\\[+0.2eM]
\hline
\end{tabular}}}
\qquad
\vspace{0.6cm}
{
\textcolor{red}{\begin{tabular}
{|p{.25cm}|p{.25cm}|p{.25cm}|
p{.25cm}|p{.25cm}|p{.25cm}|p{.25cm}|p{.25cm}|}
\hline{\bf $*$}&{\bf $*$}&{\bf $*$}&{\bf $*$}&{\bf $*$}&&&\\[+0.2eM]
\hline&&&{\bf $*$}&{\bf $*$}&{\bf $*$}&{\bf $*$}&{\bf $*$}\\[+0.2eM]
\hline{\bf $*$}&{\bf $*$}&{\bf $*$}&{\bf $*$}&{\bf $*$}&&&\\[+0.2eM]
\hline&&&{\bf $*$}&{\bf $*$}&{\bf $*$}&{\bf $*$}&{\bf $*$}\\[+0.2eM]
\hline{\bf $*$}&{\bf $*$}&{\bf $*$}&{\bf $*$}&{\bf $*$}&&&\\[+0.2eM]
\hline&&&{\bf $*$}&{\bf $*$}&{\bf $*$}&{\bf $*$}&{\bf $*$}\\[+0.2eM]
\hline{\bf $*$}&{\bf $*$}&{\bf $*$}&{\bf $*$}&{\bf $*$}&&&\\[+0.2eM]
\hline&&&{\bf $*$}&{\bf $*$}&{\bf $*$}&{\bf $*$}&{\bf $*$}\\[+0.2eM]
\hline
\end{tabular}}}
\caption{$mon_6(C_7\Box K_9)=32$ (up), $mon_7(C_8\Box K_8)=40$ (down)}
\end{center}
\label{fig1}
\end{figure}

\noindent In the following we obtain some results on the Cartesian product of complete graphs $K_n\Box K_n$.

\begin{thm}
\begin{enumerate}
\item[i)]For $n\geq 2$, $mon_{2}(K_n \Box K_n)=n$.
\item[ii)] Let for some positive integer $k$ and positive even integer $t$, $n=k\dfrac{t}{2}$, then $mon_{t}(K_n \Box K_n)=k\dfrac{t^2}{4}$.
\item[iii)] Let $k$, $t$ and $n$ be as in part (ii), then $mon_{2n-t}(K_n \Box K_n)=k(k-1)\dfrac{t^2}{4}$.
\end{enumerate}
\end{thm}

\noindent \begin{proof}
To prove $(i)$, select all of the vertices on the main diagonal of the graph as the $2$-monopoly set $M$. Every vertex of $K_n \Box K_n \setminus M$ is dominated by only $2$ monopoly vertices and each monopoly vertex contributes in activating all of its $2(n-1)$ neighbors. Hence the $2$-monopoly number of the graph is exactly $n$.

\noindent To prove  $(ii)$, let $M$ be a monopoly set for $K_n \Box K_n$ which consists all $t/2$ by $t/2$ blocks on the main diagonal. Every vertex from $K_n\Box K_n \setminus M$ has exactly $t$ neighbors in $M$ and $(ii)$ follows. For the third case, $(iii)$, let $M$ be all the vertices of $K_n \Box K_n$ except those mentioned in case $(ii)$. Each vertex of $V\setminus M$ is dominated by exactly $2n-t$ monopoly vertices. This implies the desired result.
\end{proof}

\noindent The next result concerns the line graph of regular graphs. Dynamic monopolies of this class of graphs was studied in \cite{Z}. Here we deal with the monopoly number of $L(G)$ with simple majority threshold. In other words, threshold of each vertex $v$ of $L(G)$ is considered as $\tau(v)= d_{L(G)}(v)/2$, where $d_{L(G)}(v)$ denotes the degree of $v$ in $L(G)$. Note that for each vertex $\tau(v)=k-1$, where $k$ is the regularity degree of $G$. We denote by $mo(L(G))$ the smallest cardinality of any $\tau$-monopoly in this case.

\begin{thm}
Let $G$ be a regular graph of degree $k$ with $n$ vertices. Then
\begin{center}
$ mo(L(G))\geq \dfrac{n(k-1)}{4}.$
\end{center}\label{lower}
\end{thm}

\noindent \begin{proof}
Let $M$ be a simple majority monopoly for $L(G)$. Every vertex of $M$ corresponds to an edge of $G$. Denote by $E_M$ the subset of edges in $G$ corresponding to the elements of $M$. Define the spanning subgraph $H$ of $G$ as $V(H)=V(G)$ and $E(H)=E(G)\setminus E_M$. Let $E_0\subseteq E(G)$ be any subset of edges in $G$. For any edge $e=xy$, by $d_{E_0}(e)$ we mean $d_{E_0}(x)+d_{E_0}(y)$, where $d_{E_0}(x)$ denotes the number of edges in $E_0$ incident to $x$. We have the following relations,
\begin{center}
$ \forall e=xy \in H, d_{E_M}(e)\geq k-1.$
\end{center}
Hence,
\begin{center}
$d_{E(H)}(e)=d_{E(G)}(e)-d_{E_M}(e) \leq 2k - (k-1)=k+1,$
\end{center}
which leads to
\begin{center}
$d_{E(H)}(e)=d_{H}(x)+d_{H}(y)\leq k+1.$
\end{center}
Summing up the above statement over all edges of $H$, we obtain
\begin{center}
$ \displaystyle\sum_{e=xy \in H}d_{H}(x)+d_{H}(y) \leq (k+1) \vert E(H) \vert . \qquad \quad (*)$
\end{center}
Each vertex $x$ appears $d_H(x)$ times in above sum and each time it contributes $d_H(x)$. Therefore,
\begin{center}
$\displaystyle\sum_{e=xy \in H}d_{H}(x)+d_{H}(y)=\displaystyle\sum_{x \in V(H)}d^2 (x).$
\end{center}
The right side of the equation is minimized when the graph is balanced, i.e. some vertex-degrees in it, say $\alpha$ of them equals $\lfloor \dfrac{2\vert E(H) \vert}{n} \rfloor$ and the rest, $\lceil  \dfrac{2\vert E(H) \vert}{n} \rceil$.\\
Assume that $2 \vert E(H) \vert=nq+r$, $0\leq r < n$. Then
\begin{center}
$ \displaystyle \alpha q+(n-\alpha)(q+1)=2\vert E(H) \vert$
\end{center}
which implies that
\begin{center}
$\alpha q+nq+n-\alpha q-\alpha=2\vert E(H) \vert,$
\end{center}
or
\begin{center}
$\alpha =nq+n-2\vert E(H) \vert=2\vert E(H) \vert -r+n-2\vert E(H) \vert=n-r.$
\end{center}
For simplicity, let $e=\vert E(H) \vert$. We have,
\begin{center}
$ \alpha q^2 +(n-\alpha)(q+1)^2 =(n-r)q^2 + r(q+1)^2 = (n-2e+nq)q^2 + (2e-nq)(q^2 +2q+1) = nq^2 - 2eq^2 + nq^3 + 2eq^2 + 4eq + 2e - nq^3 - 2nq^2 - nq = -nq^2 + 4eq - nq + 2e= -nq (q+1)+ 4eq + 2e \geq -n (\dfrac{2e}{n})(\dfrac{2e}{n}+1)+ 4e (\dfrac{2e}{n})+2e = -2e (\dfrac{2e}{n}+1)+8 \dfrac{e^2}{n}+2e= -4 \dfrac{e^2}{n}+8 \dfrac{e^2}{n}=4 \dfrac{e^2}{n}.$
\end{center}
In the following we demonstrate the validity of inequality in the third row of the above calculations. We have to show
$$-nq(q+1)+4eq+2e\geq -n(\frac{2e}{n})(\frac{2e}{n}+1)+4e(\frac{2e}{n})+2e.$$

\noindent We replace $q=(2e-r)/n$, simplify the inequality and obtain the following equivalent inequality.

$$-n(\frac{2e}{n}-\frac{r}{n})(\frac{2e}{n}+1-\frac{r}{n})+4e(\frac{2e}{n}-\frac{r}{n})\geq -n(\frac{2e}{n})(\frac{2e}{n}+1)+4e(\frac{2e}{n}).$$

\noindent The later inequality is equivalent to the following one.

$$\frac{2er}{n}+\frac{2er}{n}+r-\frac{r^2}{n}-4e\frac{r}{n}\geq 0.$$

\noindent The last inequality is equivalent to $nr\geq r^2$ which clearly holds. Now, from $(*)$ we get
\begin{center}
	$4 \dfrac{e^2}{n} \leq (k+1) e $
\end{center}

\noindent And this leads to
\begin{center}
$ e \leq n \dfrac{(k+1)}{4},$
\end{center}
And so
\begin{center}
$mo(L(G)) \geq \dfrac{nk}{2}- n \dfrac{(k+1)}{4}=\dfrac{n(k-1)}{4}.$
\end{center}
Note that when $k$ is odd then $n$ should be even because $kn$ is an even integer. Hence when $k$ is odd, both $n$ and $k-1$ are even. Therefore for odd $k$, $n(k-1)/4$ is an integer.
\end{proof}

\noindent Since $K_n \Box K_n = L(K_{n,n})$ then we may apply Theorem \ref{lower} for $K_n\Box K_n$ with constant threshold $n-1$. Note that when $n$ is odd we can construct a monopoly with the desired cardinality i.e. $n(n-1)/2$. In fact a circulant pattern in which each row contains exactly $(n-1)/2$ monopoly vertices, is one solution. We obtain the following result which shows that the bound of Theorem \ref{lower} is tight.

\begin{thm}
Let $n$ be an odd integer. Then $$mon_{n-1}(K_n\Box K_n) = \dfrac{n(n-1)}{2}.$$
\end{thm}

\section{Results for dynamic monopolies}

\noindent We begin this section with considering the graphs $C_n \Box K_n$ with constant threshold $\tau(v)=2$ for any vertex.

\begin{thm}
For $n\in \mathbb{N}$ we have $dyn_2(C_n\Box K_n) = \lfloor \dfrac{n}{2}+1 \rfloor$.
\end{thm}

\noindent \begin{proof}
Label the copies of the complete graphs and cycles from $K_{n}^{1}$ to $K_{n}^{n}$ and $C_{n}^{1}$ to $C_{n}^{n}$ respectively. Also, denote by $v_{ij}$ the vertices of $C_n \Box K_n$ where $i$ and $j$ refer to copies $K_{n}^{i}$ and $C_{n}^{j}$ in which $v_{ij}$ is located. If $n$ is odd, then let $D=\{v_{ii} \mid 1\leq i\leq n, i\ is \ odd \}$. $D$ is trivially a dynamic monopoly of size $\lfloor \dfrac{n}{2}+1 \rfloor$. Consequently, $dyn_2(C_n \Box K_n)\leq \lfloor \dfrac{n}{2}+1 \rfloor$. To obtain the equality, we apply the fact that from each two consecutive copies of $K_n$, one vertex is needed to construct a dynamic monopoly. That is to say, the cardinality of such a set is at least $\lfloor \dfrac{n}{2}+1 \rfloor$.\\
Now, let $n$ be an even integer. Put in $D$ the vertices
$$D=\{v_{ii} \mid 1\leq i\leq n-1, i\ is \ odd \}\cup \{v_{nn}\}.$$
We claim that $dyn_2(C_n \Box K_n)=\lfloor \dfrac{n}{2}+1 \rfloor$ in this case too. Assume to the contrary, there is a dynamic monopoly $D$ of size $\dfrac{n}{2}$ for $C_n \Box K_n$. Clearly, no two initially activated vertices in $D$ are in consecutive copies of $K_n$. Neither are they in the same $K_n$. Two cases happen for $D$:
\begin{itemize}
\item[i)] There are some vertices $v_{ir}, v_{jr} \in D$ in column $C_r$ that $j-i \geq 3$.
\\
 In this case, no vertices in the second step will be activated.
\item[ii)] There are some vertices $v_{ir}, v_{(i+2)r} \in D$ in a same column $C_r$.
\\
 In this case, we have only the vertex $v_{(i+1)r}$ activated in the second step. All other vertices will remain unactivated.
\end{itemize}
This implies that the activation process fails in either of the cases, a contradiction.
\end{proof}

\noindent We consider now $\tau(v)=3$ for all $v \in C_n \Box K_n$.

\begin{thm}
Let $n\geq3$ be any integer. Then $dyn_3(C_n \Box K_n)=n+1$.
\label{T0}
\end{thm}

\noindent \begin{proof}
Every vertex in $C_n \Box K_n$ is adjacent to at most two vertices vertically so at least one vertex in each row is required. Therefore, $\vert D \vert\geq n$. Let $D$ be a $dyn_3(C_n \Box K_n)$-set. Suppose that there are vertices that have exactly two vertical active neighbors (otherwise, the activation process will fail). These vertices become activated instantly in the second time. Others will remain unactivated unless we add some new vertex to $D$. Hence, the size of a minimum dynamic monopoly is at least $n+1$.
\\
On the other hand, let $D$ be the set of all diagonal vertices of $C_n \Box K_n$ together with an additional vertex $v_{2n}$. Then, $D$ is trivilly a dynamic monopoly of $C_n \Box K_n$ with cardinality $n+1$.
\end{proof}

\noindent The next theorem concerns general threshold assignment $t\geq 4$.

\begin{thm}
Let $m, n$ and $t$ be positive integers with $n\geq t-1$ and $t\geq 4$. Then $dyn_t(C_m \Box K_n)=m(t-2)$.
\end{thm}

\noindent \begin{proof}
For any vertex $v$ in $C_m \Box K_n$, $d(v)=(n-1)+2=n+1$ and similar to the previous argument in Theorem \ref{T0}, $v$ is adjacent to only two vertical vertices. In a certain copy of $K_n$, consider $v_{ij}$, the first vertex to be activated in its row. It has only two neighbors in $C_j$ and $n-1$ neighbors in $K_i$. Therefore, it needs at least $t-2$ initially activated vertices in the same row it is located. So, the size of a dynamic monopoly of $C_m \Box K_n$ is at least $m(t-2)$.\\
Now, we construct a dynamic monopoly for $C_m \Box K_n$. Two cases are considered:
\begin{enumerate}
\item[i)] For an even integer $m$, set
\begin{center}
$D_e =\{v_{i1}, v_{i2}, v_{i3}, \ldots, v_{i(t-2)}\mid 1 \leq i\leq m, i \ is\ odd \}$
\\
$ \bigcup \{v_{i(n-t+3)}, v_{i(n-t+4)}, \ldots, v_{i(n-1)}, v_{in}\mid 1\leq i \leq m, i \ is \ even \}.$
\end{center}
See Figure \textbf{2} for example.
\item[ii)] For an odd integer $m$, it suffices to select the following vertices to activate $C_m \Box K_n$.\\
\\
$D_o=\{v_{i1}, v_{i2}, \ldots, v_{i(t-2)}\mid 1\leq i \leq m-2, i \ is \ odd \}
 ~\cup$
 \\
 $\{v_{i(n-t+3)}, v_{i(n-t+4)}, \ldots, v_{in} \mid 1 \leq i\leq m, i \ is\ even \}~\cup~ $
 \\
 $ \{ v_{m1}, v_{m2}, \ldots, v_{m \lceil (t-2) /2 \rceil}\} ~\cup~ 
 \newline \{v_{m(n-\lfloor (t-2) /2 \rfloor)+1}, v_{m(n-\lfloor (t-2) /2 \rfloor)+2}, \ldots, v_{mn}\}. $
%\end{center}
\end{enumerate}
One can verify in an easy way that $\vert D_e \vert=\vert D_o \vert=m(t-2)$.
\end{proof}
\begin{figure}[ht]
\begin{center}
{
\textcolor{red}{\begin{tabular}
{|p{.25cm}|p{.25cm}|p{.25cm}|
p{.25cm}|p{.25cm}|p{.25cm}|p{.25cm}|p{.25cm}|p{.25cm}|p{.25cm}|}
\hline{\bf $*$}&{\bf $*$}&{\bf $*$}&&&&&&&\\[+0.2eM]
\hline&&&&&&&{\bf $*$}&{\bf $*$}&{\bf $*$}\\[+0.2eM]
\hline{\bf $*$}&{\bf $*$}&{\bf $*$}&&&&&&&\\[+0.2eM]
\hline&&&&&&&{\bf $*$}&{\bf $*$}&{\bf $*$}\\[+0.2eM]
\hline{\bf $*$}&{\bf $*$}&{\bf $*$}&&&&&&&\\[+0.2eM]
\hline&&&&&&&{\bf $*$}&{\bf $*$}&{\bf $*$}\\[+0.2eM]
\hline{\bf $*$}&{\bf $*$}&{\bf $*$}&&&&&&&\\[+0.2eM]
\hline&&&&&&&{\bf $*$}&{\bf $*$}&{\bf $*$}\\[+0.2eM]
\hline
\end{tabular}}}
\vspace{0.6cm}
\qquad
{
\textcolor{red}{\begin{tabular}
{|p{.25cm}|p{.25cm}|p{.25cm}|
p{.25cm}|p{.25cm}|p{.25cm}|p{.25cm}|p{.25cm}|}
\hline{\bf $*$}&&&&&&{\bf $*$}&{\bf $*$}\\[+0.2eM]
\hline&{\bf $*$}&{\bf $*$}&{\bf $*$}&&&&\\[+0.2eM]
\hline{\bf $*$}&&&&&&{\bf $*$}&{\bf $*$}\\[+0.2eM]
\hline&{\bf $*$}&{\bf $*$}&{\bf $*$}&&&&\\[+0.2eM]
\hline{\bf $*$}&&&&&&{\bf $*$}&{\bf $*$}\\[+0.2eM]
\hline&{\bf $*$}&{\bf $*$}&{\bf $*$}&&&&\\[+0.2eM]
\hline{\bf $*$}&&&&&&{\bf $*$}&{\bf $*$}\\[+0.2eM]
\hline&{\bf $*$}&{\bf $*$}&{\bf $*$}&&&&\\[+0.2eM]
\hline{\bf $*$}&&&&&&{\bf $*$}&{\bf $*$}\\[+0.0eM]
%\hline&{\bf $*$}&{\bf $*$}&{\bf $*$}&&&&\\[+0.0eM]
\hline
\end{tabular}}}
\caption{$dyn_5(C_{8}\Box K_{10})=24$ (up), $dyn_5(C_{9}\Box K_8)=27$ (down)}
\end{center}
\label{fig2}
\end{figure}

\noindent We shall make use of the following lemma in proving the next main result.

\noindent \begin{lemma}
Let $G$ be a graph with two threshold functions $\tau$ and $\tau'$ such that $\tau'(v)\leq \tau (v)-1$ for any vertex $v \in G$. Then $dyn_{\tau'}(G)\leq dyn_{\tau}(G)-1$.
\label{le}
\end{lemma}

\noindent \begin{proof}
Assume that $D$ is a minimum dynamic monopoly for $G$. Partition the vertices of $G$ based on the activation process runing by $D$. Let $H_0=D, H_1, \ldots, H_k$ be such a partition, i.e. $H_i$ is the set of the vertices activated at $i$th step. Hence, $\vert H_0 \vert =dyn_{\tau}(G)$. In the following we show that for each vertex $v\in G$, $H_0 \setminus v$ is a $\tau'$-dynamic monopoly for $G$. Since in $\tau'$ the thresholds is decreased by at least one and from $H_0$ we have removed exactly one vertex, then each vertex in $H_i$ is also activated at time-step $i$ with respect to the threshold function $\tau'$. Finally, all vertices except $v$, its neighbors in particular, become activated and so does $v$ itself.
\end{proof}

\begin{thm}
Let $G$ be any connected graph and $t\leq \delta(G)$ any positive integer. Let $d=dyn_t(G)$. Then
\begin{center}
$dyn_t(G \Box C_n)\leq \left \{
\begin{array}{lll}
n+d-2 \qquad \quad \quad \quad \quad \ \quad  \quad \quad t=2  & \\
nd-(n+d)+2 \quad \ \ \ \quad \ \ \   \quad \quad t=3  & \\
d+(n-2)(d-1)+d-2 \quad \quad   t\geq 4.
\end{array}
\right.$.
\end{center}
\end{thm}

\noindent \begin{proof}
Let $D$ be a $t$-dynamic monopoly for $G$. Select from the first copy the corresponding $d=dyn_t(G)$ vertices. The thresholds of the second copy has decreased to $t-1$ now. Therefore, applying Lemma \ref{le} to this copy of $G$ it can be activated by only $\vert D \vert -1$ vertices of $D$. The same happens to the third copy. It becomes active by $\vert D \vert-1$ activated vertices. There are three cases:
\begin{itemize}
\item[i)]  If $t=2$, then we need $d$ vertices for the first copy and each of the copies $G_2, G_3, \ldots, G_{n-1}$ will be activated only by one vertex. Then the last copy will be activated automatically. The number of dynamic monopoly vertices in this case therefore is $n+d-2$.
\item[ii)] If $t=3$, then $d$ vertices are needed for the first copy, $d-1$ for each of the next $n-2$ copies and only one vertex for the last copy. Hence the number of initially activated vertices in this case is
$$ d+(n-2)(d-1)+1=nd-(n+d)+2.$$
\item[iii)] If $t\geq 4$, then we need $d$ vertices for the first copy to be activated. Then each of the remained copies except the last one need $d-1$ vertices and for the last copy of $G$ we need $d-2$ vertices.
\end{itemize}
\end{proof}
\\
Adams et al. proved in \cite{ABST} that for two graphs $G$ and $H$ with a constant threshold $t$,\begin{center}
$dyn_t(G \Box H)\leq dyn_t(G)dyn_t(H).$
\end{center}
Here, we improve the bound for the Cartesian product of $K_n$ by an arbitrary graph. Also, the upper bound is strictly improved in Theorem \ref{T9} and its results.

\begin{thm}
Let $G$ be any connected graph. Then the dynamic monopoly number with threshold $t$ for $G\Box K_n$ is at most
\begin{center}
$td-(\dfrac{t^2 -3t}{2}+d).$
\end{center}
where $d$ is the $t$-dynamic monopoly number of $G$.
\end{thm}

\noindent \begin{proof}
For simplicity denote the $i th$ copy of $G$ by $G_i$. Let $D$ be a $t$-dynamic monopoly of $G$.
Consider a copy $D_1$ of the $t$-dynamic monopoly of $G$ in $G_1$, the first copy of $G$ in the Cartesian product graph. Then $G_1$ will be activated leaving the thresholds of vertices in $G_2$ to be $t-1$. Let $v$ be an arbitrary vertex of $D_1$ and $D_2=D\setminus \{v\}$. Then according to Lemma \ref{le}, $G_2$ becomes activated by selecting the corresponding vertices of $D_2$ in it. Therefore, the thresholds of $G_3$ decrease to $t-2$. Continuing so, the vertices of $G_{t-1}$ get activated and change the thresholds of $G_t$ to $1$. Any arbitrary vertex of $G_t$ now can activate this copy of $G$. Then the thresholds of all copies $G_i,\ i\geq t$ decrease to $0$ which means all of them will become activated in the next step. So, the total number of initially activated vertices is
\begin{center}
$d+d-1+d-2+\cdots +d-(t-2)+1=1+\displaystyle\sum_{i=d-(t-2)}^{d}i=1+\displaystyle\sum_{i=1}^{d}i -\displaystyle\sum_{i=1}^{d-(t-1)}i=1+\dfrac{d(d+1)}{2}-\dfrac{(d-t+1)(d-t+2)}{2}=td-(\dfrac{t^2 -3t}{2}+d).$
\end{center}
\end{proof}

\noindent But in general the bound given in \cite{ABST} is tight as shown in the following theorem.

\begin{thm}
Consider the Cartesian product of two stars, $K_{1, n}\Box K_{1, n}$, $n\geq 3$, to each vertex of which assigned a threshold $t$, $3\leq t \leq n$. Then its $t$-dynamic monopoly number is $n^2$.
\end{thm}

\noindent \begin{proof}
There are $n^2$ vertices of degree $2$ in this graph which are undoubtedly in any $t$-dynamic monopoly set. Denote by $K^1, K^2, \ldots, K^n$ each copy of $K_{1, n}$ and let
\begin{center}
$V(K^i)=\{v_{i0}, v_{i1}, \ldots, v_{in} \}, \ 1\leq i\leq n $
\end{center}
in which $v_{i0}$ is the center of $K^i$. \\
Then the set of all $v_{ij}$'s, $1\leq i, j \leq n$ will activate every $v_{i0}$, $1\leq i\leq n$. Also, each vertex $v_{0j}$, $1\leq j\leq n$ is adjacent to $n$ activated vertices and will be activated at the same time. Now, all neighbors of $v_{00}$ are activated and they are more than $n$. Hence, it will get activated too.\\
On the other hand, the dynamic monopoly number of $K_{1, n}$ with threshold $t$, $3\leq t\leq n$ is $n$. Therefore,
\begin{center}
$ dyn_t(K_{1, n}\Box K_{1, n})=n^2 =dyn_t(K_{1, n})\cdot dyn_t(K_{1, n}).$
\end{center}
\end{proof}

%\begin{figure}[!ht]
%\centerline{\includegraphics[height=4cm]{G.jpg}}
%\caption{$dyn_3(K_{1, 3}\Box K_{1, 3})=dyn_3(K_{1, 3})\cdot %dyn_3(K_{1, 3})$}
%\label{shir}
%\end{figure}
\noindent The next theorem provides an improvement of the bound of Adam et al. in \cite{ABST}.

\begin{thm}
Let $G$ and $H$ be two graphs with a constant threshold function $\tau(v)=t\geq 3$ for each one such that $dyn_t(G) \leq dyn_t(H)$. Let $D_H$ be a minimum dynamic monopoly for $H$ such that $D_H$ contains no isolated vertex. Then
\begin{center}
$dyn_t(G \Box H)\leq dyn_t(G)\cdot dyn_t(H)-\dfrac{1}{2}dyn_t(H).$
\end{center}
\label{T9}
\end{thm}

\noindent \begin{proof}
Consider the copies of $G$ in $G \Box H$ corresponding to the vertices of $D_H$. Let $P$ be a minimum dominating set in $D_H$ of cardinality say $p$. Note that $p\leq |D_H|/2$. Let $v_1$ be a vertex in $D_H$. The corresponding copy of $G$ is activated by $dyn_t(G)$ initially activated vertices. Then all neighbors of $v_1$ get threshold $t-1$ in $G \Box H$. Consider all copies of $G$ corresponding to neighbors of $v_1$ in $H[D_H]$. Each of these copies get active by $dyn_t(G)-1$ activated vertices according to Lemma \ref{le}.
Now, let $v_2$ be a next vertex in $H[D_H]$. The corresponding copy of $G$ for this vertex is activated by $dyn_t(G)$ dynamic monopoly vertices. Consider the neighbors of $v_2$ and activate the corresponding copies of $G$ by $dyn_t(G)-1$ vertices and so on. Continue this process up to $p$ vertices. Finally, the number of dynamic monopoly vertices needed in this way is
\begin{center}
$pdyn_t(G)+(dyn_t(H)-p)(dyn_t(G)-1)=dyn_t(G)dyn_t(H)-dyn_t(H)+p\leq dyn_t(G)dyn_t(H)-dyn_t(H)/2.$
\end{center}
\end{proof}

\noindent We obtain the following corollary.

\begin{cor}
If $H[D_H]$ is isomorphic to $K_{1, dyn_t(H)-1}$ for a $t \geq 3$ then
\begin{center}
$ dyn_t(G \Box H)\leq dyn_t(G) dyn_t(H)-dyn_t(H).$
\end{center}
\end{cor}

\noindent \begin{proof}
Let the central vertex of the star be $v_0$ and activated its corresponding copy of $G$ by $dyn_t(G)$ vertices. Then each copy of $G$ related to each of the neighbors of $v_0$ can be activated by $dyn_t(G)-1$ active vertices, since their thresholds have changed to $t-1$. Therefore, we have in total
\begin{center}
$dyn_t(G)+(dyn_t(H)-1)(dyn_t(G)-1)=dyn_t(G)dyn_t(H)-dyn_t(H)$
\end{center}
initially activated vertices.
\end{proof}

\begin{cor}
Assume that $dyn_t(G) \leq dyn_t(H)$. If $H[D_H]$ forms a complete graph, $K_{dyn_t(H)}$, then
\begin{center}
$dyn_t(G \Box H)\leq dyn_t(G)dyn_t(H)-\dfrac{t^2-3t+2(dyn_t(G))}{2}.$
\end{center}
\end{cor}

\noindent \begin{proof}
Consider an arbitrary vertex of $D_H$ and activate the corresponding copy of $G$ by $dyn_t(G)$ vertices. Then take another vertex in $D_H$ and activate the related copy of $G$ by $dyn_t(G)-1$ vertices. After that, consider another vertex in $D_H$. In the corresponding copy of $G$, each vertex is adjacent to two black vertices (i.e. the activated vertices until this step). Therefore, this copy will be activated by $dyn_t(G)-2$ vertices. Continue in this way and consider the $(t-1)$th vertex of $D_H$. The copy of $G$ which corresponds to this vertex can be activated by $1$ vertices since its thresholds have decreased to $2$. The last copy $G_t$ becomes activated by only one vertex of $D_H$. So, the number of initially activated vertices is
\begin{center}
$dyn_t(G)+(dyn_t(G)-1)+(dyn_t(G)-2)+\cdots +dyn_t(G)-(t-2)+1.$
\end{center}
Setting $g=dyn_t(G)$ we have
\begin{center}
$1+\displaystyle\sum_{i=g-t+2}^{g}i=1+\displaystyle\sum_{i=1}^{g}i-\displaystyle\sum_{i=1}^{g-t+1}i=1+\dfrac{g(g+1)}{2}-\dfrac{(g-t+1)(g-t+2)}{2}=gt-\dfrac{t^2 -3t+2g}{2}.$
\end{center}
\end{proof}

\noindent In the following, we consider the Cartesian product of two complete graphs. Theorem 6 in \cite{Z} states that

\begin{thm}
Let $G$ be an $r$-regular bipartite graph on $n$ vertices and $t(e)$ an assignment of thresholds to the edges of $G$. Set $t=\min\{t(e):e\in E(G)\}$. Let $D\subseteq E(G)$ be a dynamic monopoly of size $k$ in $L(G)$. Then $$k\geq \frac{n(2t-2r+2)+(2r-t)^2-4r+2t}{4} + \epsilon$$
\noindent where $\epsilon=1/4$ if $n-2r+t+1$ is an even integer and $\epsilon =0$ otherwise.\label{14}
\end{thm}

\noindent The following result extends a result in \cite{Z}.

\begin{thm}
Let $t, n$ be two fixed positive integers and $t\leq 2n-2$. Then
\begin{center}
$dyn_t(K_n \Box K_n)=\left \{
\begin{array}{ll}
\lceil \dfrac{t}{2} \rceil^2 \qquad \qquad   t\ is\ odd & \\
\dfrac{t}{2}(\dfrac{t}{2}+1) \qquad otherwise
\end{array}
\right.$
\end{center}
\label{T}
\end{thm}

\noindent \begin{proof}
First, we construct a dynamic monopoly of the desired size in each case. Let $t$ be an odd. Define
\begin{center}
$D_o =\{v_{i(n-\lfloor t/2 \rfloor +i-1)}, v_{i(n-\lfloor t/2 \rfloor +i)}, v_{i(n-\lfloor t/2 \rfloor +i+1)}, \ldots, v_{in}\mid 1\leq i \leq \lceil t/2 \rceil \}
 \cup \{v_{(n-\lfloor t/2 \rfloor +j)j}, v_{(n-\lfloor t/2 \rfloor +j+1)j}, v_{(n-\lfloor t/2 \rfloor +j+2)j}, \ldots, v_{nj}\mid 1\leq j\leq \lfloor t/2 \rfloor \}.$
\end{center}
Now assume that $t$ is an even. Let
\begin{center}
$D_e=\{v_{i(n- t/2 +i)}, v_{i(n- t/2 +i+1)}, v_{in}\mid 1\leq i \leq  t/2 \}$
\\
$ \cup \{v_{(n- t/2 +j)j}, v_{(n- t/2 +j+1)j}, v_{nj}\mid 1\leq j\leq t/2  \}.$
\end{center}
To prove the equality, let $G=K_{n,n}$ in Theorem \ref{14}.
\end{proof}

\noindent For the rest of the paper we need to generalize the above-mentioned theorem of \cite{Z}.

\begin{thm}
Let $G=(V_1, V_2)$ be a biregular bipartite graph with $\vert V_1 \vert=m, \vert V_2 \vert=n$ such that the degree of any vertex in $V_1$ (resp. $V_2$) is $r_1$ (resp. $r_2$). A fixed threshold $t \leq r_1+r_2-2$ is assigned to each edge of $G$. Let $D$ be a dynamic monopoly of size $k$ in $L(G)$. Then
$$k\geq \dfrac{mr_1 +n r_1}{2} -mn +(\dfrac{m+n+1+t-r_1 - r_2}{2})^2 - \varphi$$
where $\varphi =0$ if $m+n+1+t-r_1 - r_2$ is even and $\varphi=\dfrac{1}{4}$ otherwise.
\label{Te}
\end{thm}

\noindent \begin{proof}
The proof is similar to that of Theorem 6 of \cite{Z} in a sense that we consider an edge $e=uv$ in $H=G\setminus D$ of degree $d_H(e_1)=d_G(e_1)-d_D(e_1)\leq d_G(e_1)-t$ and eliminate its end-vertices from $H$.
Indeed, setting $H_1=H\setminus \{u,v\}$ and repeating the process, we obtain
$$|E(H_i)|\geq |E(H)| -i(r_1 + r_2-2)+i(t-1)$$
On the other hand, $|E(H_i)|\leq (m-i)(n-i)$. This leads to
%\begin{equation*}
$$|E(H)| \leq (m-i)(n-i)+i(r_1 + r_2-2)-i(t-1)$$
$$\hspace*{1.5cm}=mn+i^2 +i(r_1 +r_2 - m -n -1-t)~~(*)$$
%\end{equation*}
\noindent The latter term is minimized at $i=(m+n+1+t-r_1 - r_2)/2$. Therefore,
$$|E(H)| \leq mn -\dfrac{(m+n+1+t-r_1 - r_2)^2}{4}+\varphi.$$
Now $\varphi =0$ if $m+n+1+t-r_1 - r_2$ is even and $1/4$ otherwise.
\\
Note that $|E(G)|=(mr_1+nr_2)/2$ and $|E(H)|=|E(G)|-k$. Hence,
$$ |E(G)|-k \leq mn -\dfrac{(m+n+1+t-r_1 - r_2)^2}{4}+\varphi$$ and this implies that
$$ k\geq \dfrac{mr_1+nr_2}{2}-mn+(\dfrac{m+n+1+t-r_1 - r_2}{2})^2 -\varphi.$$
\end{proof}

\begin{cor}
Let $m, n$ and $t$ be positive integers and $t\leq m+n-2$. Then
\begin{center}
$dyn_t(K_m \Box K_n)=\left \{
\begin{array}{lll}
\lceil \dfrac{t}{2} \rceil^2 \qquad \qquad   t\ is\ odd & \\
\dfrac{t}{2}(\dfrac{t}{2}+1). \qquad otherwise.
\end{array}
\right.$
\end{center}
\end{cor}

\noindent \begin{proof}
In Theorem \ref{Te}, set $r_1 =n$ and $r_2= m$. On the other hand, taking the same vertices as those of the dynamic monopoly of $K_n \Box K_n$ in Theorem \ref{T} makes the whole graph be activated in some time steps.
\end{proof}

\begin{thm}
Let $m, n$ and $t$ be positive integers with $m<\dfrac{t}{2}$ and $n>t-m+1$. Then $dyn_t(K_m \Box K_n)=m(t-m+1)$.
\end{thm}

\noindent \begin{proof}
The proof is completely analogous to the one of Theorem 4 in \cite{Z} except that in this case the minimum value of the right hand term in the inequality $(*)$ (in proof of Theorem \ref{Te}) happens at $i=m$ since $m<t/2 < (m+n+1+t-r_1 - r_2)/2$. That is to say,
$$\hspace*{-2cm}|E(H)| \leq mn+m^2 +m(r_1 +r_2 - m -n -1-t)$$
$$=m(n+m+r_1+r_2-m-n-1-t)=m(r_1+r_2-1-t).$$
So, $$k\geq \dfrac{mr_1+nr_2}{2}-m(r_1+r_2-1-t).$$ Set now $r_1=n$ and $r_2=m$. Then,
$$k\geq mn-m(m+n-1-t)=m(t-m+1).$$
For the other side of the equality, choose $m-1, m-2, \ldots, 1$ down-side vertices in first, second, \ldots, $(m-1)$th column of $K_m\Box K_n$ respectively. Add to this set $t-(m-1), t-(m-1)-1, \ldots, 1$ right-side vertices of first, second, \ldots, $(t-m-1)$th row of the graph. Also, select $t-2(m-1)$ more vertices from the last row of the graph. Then the entire number of selected vertices in $D$ is
\begin{flushleft}
$[(m-1)+t-(m-1)]+[(m-2)+t-(m-1)-1]+[(m-3)+t-(m-1)-2]+\cdots +[1+t-(m-1)-(m-2)]+t-2(m-1)$
\\
$=t+t-2+t-4+\cdots +t-2(m-2)+t-2(m-1)$
\\
$=(m-1)t-2(1+2+\cdots+m-2)+t-2(m-1)$
\\
$=(m-1)t-(m-2)(m-1)+t-2(m-1)$
\\
$=m(t-m+1).$
\end{flushleft}
for which we counted the vertices of $i$th row and column together, where $1\leq i \leq m-1$.
\end{proof}

\section{Acknowledgment}

\noindent The authors thank the anonymous reviewers of this paper for their useful comments and suggestions.

\end{document}